\documentclass[a4paper,11pt]{article}

\usepackage[english]{babel}
\usepackage{amsmath, amssymb , latexsym, amsthm,
  epic, epsfig, rotating}

\newtheorem{definition}{Definition}[section]

\newtheorem*{thm*}{Theorem}

\newtheorem{lemma}[definition]{Lemma}
\newtheorem*{lemma*}{Lemma}

\newcommand{\bbZ}{\mathbb{Z}}
\newcommand{\bbR}{\mathbb{R}}

\renewcommand{\phi}{\varphi}
\renewcommand{\epsilon}{\varepsilon}

\begin{document}
\sloppy

\title{There is no minimal action of $\bbZ^2$ on the plane}
\author{Fr\'ed\'eric Le Roux\footnote{Laboratoire de math\'ematique (CNRS UMR 8628),  Universit\'e Paris Sud,
91405 Orsay Cedex, France.}
}

\maketitle


\normalsize




\bigskip
\bigskip

\section{Introduction}

A group acting on a topological space is said to act \emph{minimally} if every orbit is dense.
Let $g_{1} : (x,y) \mapsto (x+1,y)$ be the horizontal translation of the plane, $g_{2} : (x,y) \mapsto (x,y+1)$ the vertical translation, and $g_{3} : (x,y) \mapsto (x,y) + \vec v$ a translation by a vector $\vec v$ with irrational coordinates and irrational slope. 
 Then the $\bbZ^3$ action generated by $g_{1},g_{2},g_{3}$  on the plane is minimal. On the other hand there is no minimal homeomorphism of the plane: this is a  consequence of Brouwer plane translation theorem.  In other words the group $\bbZ$ does not act minimally. Thus it is a natural question to determine whether the group $\bbZ^2$ admits a minimal plane action. In this paper, we answer this question negatively.

\begin{thm*}
There is no minimal action of $\bbZ^2$ by homeomorphisms of the plane.
\end{thm*}

Here is another motivation. Assume $\bbZ^2$ acts on the plane, let $g_{1},g_{2}$ be generators of the action. If $g_{1}$ is conjugate to a translation, then the quotient $\bbR^2/g_{1}$ is homeomorphic to the infinite annulus, and $g_{2}$ induces an annulus homeomorphism. If the $\bbZ^2$ action was minimal, then so would be the induced annulus homeomorphism. But this would contradict Le Calvez-Yoccoz's theorem.

\begin{thm*}[Le Calvez-Yoccoz, \cite{LeCalvezYoccoz97}]
There is no minimal homeomorphism of the infinite annulus.
\end{thm*}

Thus the non-existence of minimal $\bbZ^2$ action may be seen as a (slight) generalisation of Le Calvez-Yoccoz's theorem.

In contrast, there exist actions of $\bbZ$ on the plane with some dense orbits; however, due again to Brouwer's theorem, such an action cannot be free.
There also exist free transitive actions of $\bbZ^2$. Both examples may be constructed using transitive skew products on the infinite annulus
 (see~\cite{Besicovitch51}). I do not know if there exists a transitive free action of $\bbZ^2$ for which  no element is conjugate to a translation.

\subsection*{Sketch of the proof}

We argue by contradiction. Let $g_{1}, g_{2}$ be two commuting homeomorphisms of the plane that generates a minimal action of $\bbZ^2$. For simplicity, let us assume that these homeomorphisms preserve the orientation. We first prove (Lemma~\ref{lem.minimal-free}) that the action of $\bbZ^2$ is free; this is essentially due to the non-existence of minimal actions of $\bbZ$. Applying Le Calvez-Yoccoz's theorem, we see that no element of the action is conjugate to a translation. Let $D$ be a small disk. Since $g_{1}g_{2}$ is not conjugate to a translation, the sequence of positive iterates of $D$ under $g_{1}g_{2}$ accumulates somewhere in the plane. The same holds for $g_{1}g_{2}^{-1}$.  From this, using the minimality of the action, we are able to construct a ``periodic pseudo-orbit'' for $g_{1}$ (more precisely, a \emph{periodic disk chain}). This is too much recurrence for a fixed point free orientation preserving homeomorphism of the plane, as encapsuled by John Franks's \emph{disk chain lemma}, and we get a contradiction.
\emph{}

\section{Background material}

\subsection*{Brouwer homeomorphisms}
A \emph{Brouwer homeomorphism} is a fixed point free, orientation preserving homeomorphism of the plane. Let $h$ be a Brouwer homeomorphism. A  disk is \emph{free} for $h$ if it does not meet its image. The \emph{free disk lemma} says that every free disk $D$ is disjoint from all its iterates $h^n(D), n \neq 0$. In particular, a Brouwer homeomorphism has no periodic orbit, and no dense orbit. Another immediate consequence is that there is no minimal orientation preserving homeomorphism of the plane. See for example~\cite{Brouwer12, HommaTerasaka53,Franks92,Guillou94,BeguinLeRoux03} for more details about Brouwer homeomorphisms.

A  couple $(x,y)$ of points in the plane is said to be \emph{singular} for the Brouwer homeomorphism $h$ if there exist a sequence $(z_{k})$ of points converging to $x$, and a sequence $(n_{k})$ of integers tending to $+\infty$, such that the sequence of points $(h^{n_{k}}(z_{k}))$ converges to $y$.
The \emph{singular set} of $h$ is the set of singular couples. It has the following properties. The singular set of $h^n$ is equal to the singular set of $h$ for every $n\neq 0$. Furthermore, this set is empty if and only if $h$ is conjugate to a translation. From this we see that if some power $h^n, n \neq 0$ is conjugate to a translation, then so is $h$. In particular, concerning our problem, Le Calvez-Yoccoz's theorem may be reformulated as follows.

\begin{thm*}[Le Calvez-Yoccoz]
There is no minimal action of $\bbZ^2$ on the plane for which some element is conjugate to a translation.
\end{thm*}

The \emph{positive singular set} of $h$ is the first projection of the singular set, that is, the set of points $x$ such that $(x,y)$ is a singular couple for some point $y$.

A \emph{disk chain} for $h$ is a sequence  $(D_i)_{i=1..k}$ of pairwise disjoint free disks such that
for every $i=1, \cdots,  k-1$, some positive iterate of the disk $D_{i}$ meets $D_{i+1}$.
The disk chain is \emph{periodic} if some positive iterate of $D_{k}$ meets $D_{1}$.
The \emph{disk chain lemma} is a  generalization of the free disk lemma:
\begin{lemma*}[John Franks,~\cite{Franks92}]
A Brouwer homeomorphism admits no periodic disk chain.
 \end{lemma*}

In all our statements about disks we will assume that the disks are open.
This does not really matter here (in particular, all the previous statements also hold for closed disks).
Note that the disk chain lemma is still true if one replaces euclidean disks by \emph{topological disks}, that is, simply connected open subsets of the plane.

\subsection*{Orientation-reversing Brouwer homeomorphisms}

In this paper we will consider actions by homeomorphisms of the plane, and we will not assume that these homeomorphisms preserve the orientation. Thus we also need to collect some results about orientation reversing homeomorphisms. We begin by a lemma that was certainly known to  K\'er\'ekj\'art\'o.

\begin{lemma}
\label{lem.racine-translation}
Let $g$ be a homeomorphism of the plane such that some power $g^p$ is a translation. Then $g$ is conjugate either to a translation or to the map $(x,y) \mapsto (x+1, -y)$.
\end{lemma}
\begin{proof}
Under the assumptions it is easy to see that the map $g$ acts properly discontinuously on the plane, and thus the quotient space $\bbR^2/g$ is a (Hausdorff) surface. Moreover, this surface has a $p$-fold covering by the annulus $\bbR^2/g^p$. The classification of surfaces now tells us that $\bbR^2/g$ is homeomorphic either to the infinite annulus $\bbR^2/(x,y) \mapsto (x+1,y)$,  or to the  Mobius band $\bbR^2/(x,y) \mapsto (x+1,-y)$.  This provides the two cases, since $g$ is an automorphism of the covering map   $\bbR^2/g^p \to \bbR^2/g$.
\end{proof}

Fortunately, Brouwer theory has been adapted to the orientation-reversing case by M. Bonino. We will make use of  the following \emph{orientation-reversing free disk lemma} (see~\cite[Lemma 5.2]{Bonino04}).
Let $g$ be an orientation reversing homeomorphism of the plane with no periodic point of period $2$,
 and $D$ be a  disk which is disjoint from both $g(D)$ and $g^2(D)$.
Then $D$ is disjoint from all its iterates $g^n(D), n\neq0$. In particular, there is no orientation reversing minimal homeomorphism of the plane.

Next we notice that the argument given by Le Calvez and Yoccoz on the annulus works equally well on the Mobius band.
\begin{thm*}[essentially due to Le Calvez-Yoccoz]
There is no minimal homeomorphism of the Mobius band.
\end{thm*}

\begin{proof}
Assume there is some minimal homeomorphism $h$ of the Mobius band $M$.  The homeomorphism $h$ extends to a homeomorphism of  the one-point compactification $M \cup \{\infty\}$ that fixes the point $\infty$. If $h$ is minimal then this point is not a source nor a sink, and no neighbourhood of it can contain a whole orbit (see the precise hypotheses in \cite{LeCalvezYoccoz97}). Then Le Calvez-Yoccoz's index theorem applies, and tells us that for some positive 
 iterate of  $h$, the fixed point index of $\infty$ is non-positive. On the other hand, this iterate is isotopic to the identity (see~\cite{Hamstrom65}) and we may apply Lefchetsz formula:  the sum of its fixed points indices is equal to $1$ (the Euler caracteristic of $M \cup \{\infty\}$, topologically a projective plane). This shows that $h$ has some periodic orbit on $M$, in contradiction with the hypothesis that it is minimal.
\end{proof}

\subsection*{Minimal actions are free}
Let $G$ be a group of homeomorphisms of the plane. Recall that $G$  is \emph{free} if every non trivial element of $G$ is fixed point free.
The following lemma explains why Brouwer theory is useful for minimal actions.
 \begin{lemma}\label{lem.minimal-free}
Assume $G$ is isomorphic to $\bbZ^2$. If the action of $G$ is minimal then it is free.
\end{lemma}

An element $g_{1}$ of $G$ will be called \emph{primitive} if there is no element $g' \in G$ and integer $p > 1$ such that $g_{1} = g'^p$. Equivalently, when $G$ is isomorphic to $\bbZ^2$, $g_{1}$ is primitive if there exists $g_{2} \in G$ such that $\{g_{1},g_{2}\}$ generates $G$.

\begin{proof}
Assume, by contradiction, that some element of $G$ has a fixed point $x_{0}$. This element is equal to  $g_{1}^p$ for some primitive element $g_{1}$.  Since $G$ is abelian, every point of the $G$-orbit of $x_{0}$ is fixed by $g_{1}^p$. If the action is minimal then the $G$-orbit of $x_{1}$ is dense, and thus $g_{1}^p$ is the identity.
By a theorem of K\'er\'ekj\'art\'o (see~\cite{ConstantinKolev94}), $g_{1}$ is conjugate to an isometry.  Using the classification of plane isometries we see that $g_{1}$ itself must have a fixed point. Again, since the group is abelian and acts minimally, we get that $g_{1}$ is the identity.
  Now the group $G/g_{1}$ is isomorphic to $\bbZ$ and acts minimally on the plane, which contradicts the (orientation preserving or reversing) free disk lemma.
\end{proof}

\section{The orientation-preserving case}

\begin{proof}[Proof of the theorem in the orientation-preserving case]
Let $G$ be a group of orientation preserving homeomorphisms of the plane.
We argue by contradiction, assuming that $G$ is isomorphic to $\bbZ^2$ and acts minimally on the plane.
According to lemma~\ref{lem.minimal-free}, the non-trivial elements of $G$ are Brouwer homeomorphisms.
According to Le Calvez-Yoccoz's theorem, they are not conjugate to a translation.
Let $g_{1}, g_{2}$ be a set of generators of $G$. The Brouwer homeomorphism $g_{1}g_{2}$ is not conjugate to a translation, thus it has a non-empty positive singular set. Since the action is minimal, and the singular set is a conjugacy invariant, the positive singular set of $g_{1}g_{2}$ is dense in the plane\footnote{Note that there is no contradiction here: examples of Brouwer homeomorphisms with a dense positive singular set are described in~\cite{HommaTerasaka53} or~\cite[page 18]{LeRoux04}.}. Choose any sufficiently small (open) disk $D$, so that $D$ is disjoint from its image under each of the four elements $g_{1}, g_{2}, g_{1}g_{2}, g_{1}g_{2}^{-1}$ (remember that the action is free). Let $(x,y)$ be a singular couple, for $g_{1}g_{2}$, with $x$ in  $D$. By minimality of the action there exists some $g = g_{1}^{n_{1}} g_{2}^{n_{2}}$ such that $g(y)$ belongs to $D$, or equivalently $y$ belongs to  $g^{-1}(D)$. Since $(x,y)$ is singular for $g_{1}g_{2}$, there exists some arbitrarily large positive integer $N$ such that $(g_{1}g_{2})^N (D)$ meets $g^{-1}(D)$, that is,  the topological disk 
$$g_{1}^{N+n_{1}}g_{2}^{N+n_{2}} (D)$$
 meets $D$. 
In particular we have found some positive integers $N_{1}, N_{2}$ such that $g_{1}^{N_{1}}g_{2}^{N_{2}}(D)$ meets $D$.

Arguing similarly with $g_{1}g_{2}^{-1}$, we find another couple of positive integers, say  $N'_{1}, N'_{2}$, such that $g_{1}^{N'_{1}}g_{2}^{-N'_{2}}(D)$ meets $D$.

We will now get a contradiction by finding some $g_{2}$-iterates of $D$ that constitute a disk chain for the Brouwer homeomorphism $g_{1}$. More precisely we get the conclusion out of the following lemma.

\begin{lemma}\label{lem.return-time}
Let $g_{1}, g_{2}$ be two commuting Brouwer homeomorphisms. For any  disk $D$, consider the \emph{return-time set} 
$$
R(D,g_{1},g_{2}) = \{(n_{1}, n_{2}) \in \bbZ^2 ,     g_{1}^{n_{1}} g_{2}^{n_{2}} (D) \cap D \neq \emptyset         \}.
$$
Assume this set does not contain  $(1,0)$ nor $(0,1)$: in other words the disk $D$ does not meet its image under $g_{1}$ and $g_{2}$. Then 
$R(D,g_{1},g_{2})$ is either contained in the set 
$$E_{+} = \{(n_{1}, n_{2}), \ \ n_{1}n_{2} >0  \}$$
 or in the set 
$$E_{-} = \{(n_{1}, n_{2}), \ \  n_{1}n_{2} <0  \}.$$
\end{lemma}
See~\cite{BargeFranks93} for some nice results about the set $R(D,g)$ in the case of a free $\bbZ$-action.

\begin{proof}[Proof of the lemma]
Under the assumptions of the lemma, first note that the set $R(D,g_{1},g_{2})$ does not contain any couple of the form $(0,n)$ nor $(n,0)$ with $n \neq 0$: indeed  this is the content of the free disk lemma when applied to $g_{1}$ and $g_{2}$.

We argue by contradiction, assuming that the set $R(D,g_{1},g_{2})$ contains $(N_{1},N_{2})$ and $(N'_{1}, -N'_{2})$ where $N_{1}, N_{2}, N'_{1}, N'_{2}$ are positive integers (as was the case in the course of the proof of the theorem). 
We write
$$N_{2} = n_{2}p , \ \ \  N'_{2} = n'_{2} p $$
where $p$ is the greatest commun divisor of $N_{2}, N'_{2}$. 
Consider the (cyclic) sequence of topological disks $(D_{i})_{i \in \bbZ/(n_{2}+n'_{2}+1)\bbZ}$ consisting in the following $g_{2}$-iterates of $D$:
$$  D , \ g_{2}^{N'_{2}}(D), \ g_{2}^{2N'_{2}}(D), \ \dots,  \ g_{2}^{n_{2} N'_{2}}(D) = g_{2}^{n'_{2} N_{2}}(D), \ g_{2}^{(n'_{2}-1) N_{2}}(D), \ \dots , \ g_{2}^{N_{2}}(D).$$

Since $p n_{2}n'_{2} = n'_{2} N_{2} = n_{2} N'_{2}$ is  the least commun multiple of $N_{2}$ and $N'_{2}$, all the powers of $g_{2}$ involved in the sequence are distinct. Thus the disks in the sequence are pairwise disjoint, otherwise the set $R(D,g_{1},g_{2})$ would contain some couple of the form $(0,n)$. By assumption 
 $R(D,g_{1},g_{2})$ contains $(N_{1},N_{2})$ and $(N'_{1}, -N'_{2})$. This means that 
for each disk $D_{i}$ in this cyclic sequence,  the image of $D_{i}$ under either $g_{1}^{N'_{1}}$ or $g_{1}^{N_{1}}$ meets the next disk $D_{i+1}$. Thus this sequence constitutes a periodic disk chain for the Brouwer homeomorphism $g_{1}$, in contradiction with Franks's lemma.
\end{proof}

This completes the proof of the theorem in the orientation-preserving case.
\end{proof}

\section{The orientation-reversing case}
We will explain how to  modify the arguments of the previous section.

\begin{proof}[Proof of the theorem in the orientation-reversing case]
Assume as before that $G$ is a group of homeomorphisms of the plane, that is isomorphic to $\bbZ^2$, and acts minimally. Also assume that it contains some orientation reversing elements. Lemma~\ref{lem.minimal-free} applies, so again the elements of $G$ are fixed point free.
No element of $G$ is conjugate to a translation: otherwise $G$ would contain some primitive element conjugate either to a translation or to the map $(x,y) \mapsto (x+1, -y)$ (see Lemma~\ref{lem.racine-translation}), we would get a minimal homeomorphism of the annulus or the Mobius band, both situations contradicting Le Calvez-Yoccoz's results.

 The set $G^+$ of orientation-preserving elements of $G$ is a subgroup of index $2$ which is again isomorphic to $\bbZ^2$.
We now consider a basis $\{g_{1}, g_{2}\}$ of $G$ such that $G^+$ contains $g_{1}$ (and not $g_{2}$).
Let $D$ be a disk which is disjoint from its images under the five maps $g_{1}, g_{2}, g_{2}^2, g_{1}g_{2}^2, g_{1}g_{2}^{-2}$.
The last two maps are (orientation preserving) Brouwer homeomorphisms whose positive singular set is dense: as above we find some
 positive integers $N_{1}, N_{2}, N'_{1}, N'_{2}$ such that the return-time set $R(D,g_{1},g_{2})$ (defined as in Lemma~\ref{lem.return-time}) contains $(N_{1},N_{2})$ and $(N'_{1}, -N'_{2})$.

Since $D, g_{2}(D), g_{2}^2(D)$ are mutually disjoint, we may apply Bonino's orientation reversing free disk lemma: we get that $D$ is disjoint from all its $g_{2}$-iterates. Now, as in Lemma~\ref{lem.return-time}, we can make a periodic disk chain for $g_{1}$ out of the $g_{2}$-iterates of $D$. Since $g_{1}$ is a Brouwer homeomorphism, this again contradicts Franks's Lemma.
\end{proof}

\end{document}